\documentclass[10pt,a4paper]{article}
\usepackage[margin=1in]{geometry}

\usepackage{amscd}
\usepackage{amsfonts}
\usepackage{amsmath}
\usepackage{amssymb}
\usepackage{amstext}
\usepackage{amsthm}
\usepackage{color}
\usepackage[nodate]{datetime}
\usepackage{dsfont}
\usepackage{enumerate}
\usepackage{fancyhdr}
\usepackage{framed}
\usepackage{graphicx}
\usepackage{hyperref}
\usepackage[utf8]{inputenc}
\usepackage{mathrsfs}
\usepackage{multicol}
\usepackage{pifont}
\usepackage{psfrag}
\usepackage{stmaryrd}
\usepackage[normalem]{ulem}
\usepackage{url}
\usepackage{verbatim}
\usepackage{wasysym}

\makeatletter
\renewcommand{\paragraph}{%
  \@startsection{paragraph}{4}%
  {\z@}{0.5ex \@plus 1ex \@minus .2ex}{-1em}%
  {\normalfont\normalsize\bfseries\maybe@addperiod}%
}
\newcommand{\maybe@addperiod}[1]{#1\@addpunct{.}}
\makeatother

\newcommand{\E}{\mathbb{E}}

\newcommand{\I}{\mathds{1}}

\newcommand{\N}{\mathbb{N}}

\newcommand{\PP}{\mathbb{P}}

\newcommand{\Z}{\mathbb{Z}}

\newcommand{\Ir}{\mathcal{I}}
\newcommand{\sleep}{\mathfrak{s}}
\newcommand{\dd}{\mathrm{d}}

\theoremstyle{plain}
\newtheorem{theorem}{Theorem}

\newtheorem{lemma}[theorem]{Lemma}

\newtheorem*{assumption*}{Assumption}
\newtheorem*{lemma*}{Lemma}
\newtheorem*{proposition*}{Proposition}
\newtheorem*{theorem*}{Theorem}
\theoremstyle{definition}
\newtheorem*{definition*}{Definition}
\theoremstyle{remark}

\newtheorem*{remark*}{Remark}

\renewcommand{\geq}{\geqslant}
\renewcommand{\ge}{\geqslant}

\renewcommand{\le}{\leqslant}

\hypersetup{colorlinks=true,urlcolor=blue,linkcolor=black,citecolor=black}
\newcommand*{\doi}[1]{\href{http://dx.doi.org/\detokenize{#1}}{doi}}

\usepackage{setspace}
\begin{document}

\title{Universality and Sharpness in Activated Random Walks}

\author{
Leonardo T. Rolla,
Vladas Sidoravicius,
Olivier Zindy
\\
\small
University of Buenos Aires, NYU-Shanghai, University Paris Pierre et Marie Curie
}

\maketitle

\begin{abstract}
We consider the Activated Random Walk model in any dimension with any sleep rate and jump distribution and ergodic initial state.
We show that the stabilization properties depend only on the average density of particles, regardless of how they are initially located on the lattice.
\end{abstract}

\section{Introduction}

The concept of self-organized criticality was introduced in the late 80's to explain the emergence of critical behavior in steady states without fine tuning of system parameters~\cite{BakTangWiesenfeld88}. Intrinsic relationships between this phenomenon and that of ordinary phase transition started to unfold in the late 90's with a new paradigm: the self-organized critical behavior of a driven-dissipative system is related to ordinary criticality of a corresponding fixed-energy system that uses the same relaxation mechanism~\cite{DickmanMunozVespignaniZapperi00}. A central issue is the \emph{density conjecture}: the typical density $\zeta_s$ in the steady state of the driven-dissipative system arguably coincides with the threshold density $\zeta_c$ of the fixed-energy system. Ten years later, the density conjecture was shown to be false~\cite{FeyLevineWilson10a} for the Abelian Sandpile Model (ASM), which has been attributed to the fact that the ASM dynamics does not wipe out details of the initial condition by the time a configuration becomes explosive~\cite{JoJeong10}, see below for more details.

Stochastic Sandpiles (SSM) and Activated Random Walks (ARW) were introduced as alternative non-deterministic relaxation mechanisms.
The long-range space-time correlations caused by conservation of particles and the lack of an algebraic structure similar to the ASM make the mathematical analysis of these models very challenging.
It took two decades for the first rigorous results regarding stability properties of these systems to appear in the literature~\cite{DickmanRollaSidoravicius10,RollaSidoravicius12}. Considerable progress for the ARW has been made in the past three years~\cite{BasuGangulyHoffman18,CabezasRollaSidoravicius14,CabezasRollaSidoravicius18,RollaTournier18,SidoraviciusTeixeira17,StaufferTaggi18,Taggi16}, with the introduction of a number of \textsl{ad hoc} techniques and tools. Some of these tools were sensitive to assumptions about the initial state, requiring independence, light tails, etc.

In this paper we study the ARW and prove that the critical density $\zeta_c$ is well-defined and separates two entire families of spatially ergodic states: those whose density $\zeta$ is below $\zeta_c$, for which configurations are a.s.\ stabilizable, and those having density above $\zeta_c$, for which configurations are a.s.\ explosive.
One consequence is that special assumptions about the initial state that were needed in many of the recently-introduced mathematical techniques can now be waived.
Another consequence is to support the general belief that the ARW has much better mixing properties than the ASM. The result and technique introduced in this paper do not yield a proof of the density conjecture for the ARW, but may give one step in that direction, as discussed further below.
 

\paragraph{Sharpness and self-organized criticality}

Models of avalanches became a standard example of self-organized criticality in the context of non-equilibrium steady states.
Unlike usual statistical mechanics systems, these models are not explicitly equipped with a tuning parameter at which a phase transition is observed.
Instead, they are expected to spontaneously drive themselves to a critical steady state, featuring characteristics of critical systems such as power law statistics and scale invariance.

A common setup involves a reaction-diffusion evolution, where sites of a graph contain particles which dissipate according to certain rules, until stabilizing.
Three different relaxation procedures have been used in the study of the above phenomenon: in the deterministic sandpile model, sites with at least $2d$ particles send one to each neighbor; in the stochastic sandpile, sites with at least $2$ particles send $2$ particles to neighbors chosen at random; in the activated random walk model, sites with active particles send one particle to a neighbor chosen at random, and particles can become passive with probability described by a parameter $\lambda$ if they are alone.
All these systems contain mechanisms that cause both spread of activity and a tendency of this activity to die out, and the system behavior is determined by the balance between these two factors.

Self-organized criticality appears in the corresponding driven-dissipative dynamics: particles are added to the bulk of a large finite box, and absorbed at its boundary during relaxation, following one of the above-mentioned mechanisms. A particle is added only after the system globally stabilizes. In this dynamics, when the average density $\zeta$ inside the box is too small, mass tends to accumulate. When it is too large, there is intense activity and substantial dissipation at the boundary. Within this setting, the system is attracted to a critical state with an average density $\zeta_s$.

A new paradigm was introduced in~\cite{DickmanMunozVespignaniZapperi00}, arguing that self-organized criticality in these systems is related to ordinary phase transition. More precisely, the corresponding conservative systems in infinite volume, where the density $\zeta$ is kept constant, exhibit ordinary phase transition and their critical behavior is closely related to properties of the self-organized critical system described above. In particular, there is a threshold density $\zeta_c$ such that the infinite-volume dynamics should fixate for $\zeta<\zeta_c$ and remain active for $\zeta>\zeta_c$, and moreover $\zeta_c$ should coincide with the driven-dissipative stationary density $\zeta_s$.
Since then, a rich literature appeared, exploring this relation and the principles behind it.

Later on, rigorous results and precise large-scale simulation showed that the \emph{density conjecture} $\zeta_c=\zeta_s$ is false in general, at least for the deterministic sandpile model~\cite{FeyLevineWilson10a,FeyMeester15}, indicating that the relation between driven-dissipative and conservative systems is much more subtle.
This was attributed to the fact that the ASM is very sensitive to the initial state~\cite{JoJeong10}.
Indeed, for every $d<\zeta<2d-1$ there are spatially ergodic states with average particle density $\zeta$ that are explosive and some others with same density which are stabilizable~\cite{Fey-denBoerRedig05}, see \S2 for definitions.


This discovery had two implications.
First, it increased interest in the mathematical properties of the ASM and its intricate behavior, and propelled debate about how to recover the density conjecture for this system~\cite{Levine15,PoghosyanPoghosyanPriezzhevRuelle11}.
Second, it increased physical interest in other models such as SSM and ARW, which are supposed to have stronger mixing properties, for the study of self-organized criticality and analysis of avalanche statistics.
We show that $\zeta_c$ is well-defined for the ARW, providing some support for the latter claim.
We now state the main result of this paper, postponing precise definitions to \S\ref{sec:definitions}.

\begin{theorem}
\label{thm:main}
Consider the Activated Random Walk model on the usual graph $\Z^d$ for fixed $d\ge 1$ with given sleep rate $\lambda$, and given jump distribution $p(\cdot)$.
Assume the support of $p(\cdot)$ generates $\Z^d$ and not a sublattice.
There is a number $\zeta_c$ such that, for any spatially ergodic distribution $\nu$ supported on active configurations with average density $\zeta$, a configuration sampled from $\nu$ is a.s.\ stabilizable if $\zeta<\zeta_c$ and a.s.\ explosive if $\zeta>\zeta_c$.
\end{theorem}


\paragraph{Heuristics for the density conjecture}

Consider the driven-dissipative system corresponding to the ARW and starting with the empty configuration on a large box $\Lambda$.
Then add a number $u \, \cdot|\Lambda|$ of particles at random locations, one by one, letting the system stabilize in between.
By the Abelian property, this is the same as starting from an i.i.d.\ $\mathrm{Poisson}(u)$ number of particles at every site and then letting the configuration stabilize.
After stabilization, some particles will exit the box, and a number $\zeta(u)\cdot|\Lambda|$ will be retained.
Since every ergodic state with density $\zeta<\zeta_c$ is stabilizable, for such densities the driven-dissipative system should not be losing a macroscopic amount of mass during stabilization, and $\zeta$ should increase at rate $\frac{\dd \zeta}{\dd u} \approx 1$, at least until $\zeta\approx\zeta_c$.
On the other hand, if every state with density $\zeta>\zeta_c$ is explosive, then for $\zeta>\zeta_c$ the system should lose all the excess before $u$ increases macroscopically, down to $\zeta\approx\zeta_c$.
Combining these two facts, the graph of $\zeta \times u$ should 
increase at rate~1 for $u\in[0,\zeta_c]$ and remain constant for $u\in[\zeta_c,\infty)$.
That is, $\zeta$ should be given by $\zeta(u)=\min\{u,\zeta_c\}$, which in particular implies the density conjecture $\zeta_s := \lim\limits_{u\to\infty} \zeta(u) = \zeta_c$.

The converse argument helps to explain how failure of the ASM to satisfy the density conjecture may be related to the lack of a sharp transition at a unique point $\zeta_c$.
Consider the same setup as above but using the ASM as relaxation mechanism.
Again, adding $u \, \cdot|\Lambda|$ slowly to an empty configuration is the same as starting from an i.i.d.\ $\mathrm{Poisson}(u)$ configuration and stabilizing after.
Let $\zeta_c$ denote the threshold density for i.i.d.\ Poisson distributions, so $\zeta_c$ is the value of $u$ at which the dissipation mechanism starts to throw a positive proportion of the particles out of the box $\Lambda$.
Fig.~3 in~\cite{FeyLevineWilson10a} is a sample plot of $\zeta \times u$, showing that the curve is close to $\zeta(u)=\min\{u,\zeta_s\}$ but not quite.
In fact, it is off by a small but macroscopic difference, which went unnoticed in the physics community for many years.
In graphs such as $\Z^2$, the curve $\zeta \times u$ rises at rate~$1$ until $\zeta_c>\zeta_s$, and only then it starts to decrease smoothly towards $\zeta_s$.
This is only possible because there are states with density $\zeta \approx \zeta_c > \zeta_s$ which are still stabilizable.


\paragraph{Mathematical overview}


Several non-trivial bounds for $\zeta_c$ were proved in the past three years.
For $d=1$, it was proved in~\cite{RollaSidoravicius12} that $\zeta_c>0$ for all $\lambda$ and $\zeta_c \to 1$ as $\lambda\to\infty$.
For $d \ge 2$ and $\lambda=\infty$, it was also shown in~\cite{Shellef10} that $\zeta_c>0$ and in~\cite{CabezasRollaSidoravicius14,CabezasRollaSidoravicius18} that $\zeta_c \ge 1$.
For $d \ge 2$ and $\lambda>0$ it was shown in~\cite{SidoraviciusTeixeira17} that $\zeta_c>0$, assuming short-range unbiased jump distributions.
This was generalized to general jump distributions in~\cite{StaufferTaggi18}, where it was also shown that $\zeta_c \to 1$ as $\lambda\to\infty$.
It was proved in~\cite{AmirGurel-Gurevich10,Shellef10} that $\zeta_c \le 1$ in any dimension for any $\lambda$.
For biased jump distributions, it was shown in~\cite{Taggi16} that, on $d=1$, $\zeta_c<1$ for every $\lambda$ and $\zeta_c\to 0$ as $\lambda \to 0$, and on $d \ge 2$ that $\zeta_c<1$ for small $\lambda$.
The picture on $d \geq 2$ was extended in~\cite{RollaTournier18} by showing show that $\zeta_c<1$ for every $\lambda$ and $\zeta_c \to 0$ as $\lambda\to 0$.
For unbiased jumps, it was shown that $\zeta_c\to 0$ as $\lambda \to 0$, in~\cite{BasuGangulyHoffman18} for $d=1$ and~\cite{StaufferTaggi18} for $d \ge 3$.
See~\cite{Rolla15} for a detailed account.

While many of the proofs were robust with respect to the initial state, others were sensitive to it.
Some results in \cite{AmirGurel-Gurevich10,CabezasRollaSidoravicius14,CabezasRollaSidoravicius18,StaufferTaggi18,Taggi16} required an i.i.d.\ field, the proofs in~\cite{SidoraviciusTeixeira17} are presented for i.i.d.\ Poisson (but may be adapted for states with some spatial mixing and light tails), and some of the results in~\cite{StaufferTaggi18} and~\cite{Taggi16} required the initial state to be i.i.d. Bernoulli.
With Theorem~\ref{thm:main}, all these conditions can be waived, and arguments having special requests for the initial state now produce lower and upper bounds valid for any other ergodic initial state.
A good example is the case of biased walks in dimension $d \geq 2$.
The proof in~\cite{Taggi16} that $\zeta_c<1$ for all $\lambda>0$ and i.i.d.\ Bernoulli initial state was extended to an arbitrary i.i.d.\ field in~\cite{RollaTournier18}, however such an extension now follows directly from this general property.

%

\section{Definitions and tools}
\label{sec:definitions}

In this section, we recall usual tools for this system and fix the notation. We quickly describe a continuous-time evolution and the event of fixation, the site-wise representation and the notion of stabilization with its main properties, the relation between stabilization and fixation, and finally we state the mass-transport principle.

To be consistent with most of the existing literature, we describe the ARW evolution as a continuous-time stochastic process.
The ARW system starts with active particles placed in $\Z^d$ according to a distribution $\nu$, and evolves as follows.
Active particles perform independent continuous-time random walks on $\Z^d$ with translation-invariant jump distribution $p(x,y)=p(y-x)$, and switch to passive state at rate $\lambda>0$ when they are alone on a site.
Passive particles do not move, and are reactivated immediately when visited by another particle.
The law of this evolution is denoted $\PP^\nu$.
By fixation we mean that the dynamics eventually halts at any finite region, and non-fixation is the opposite event.

In the sequel we introduce notation for this continuous-time evolution, and briefly describe the toppling operators and their properties, referring the reader to~\cite{Rolla15} for detailed explanations.
One important difference in notation between here and~\cite{Rolla15} is that here we do not necessarily start from a zeroed odometer.
One difference between the terminology of~\cite{Rolla15} and that of~\cite{RollaSidoravicius12,SidoraviciusTeixeira17} is that we do not replace sleep instructions by neutral ones, but instead allow for the toppling of a site containing a sleepy particle.

\paragraph{Notation and continuous-time evolution}

Let $\N_0 = \{0,1,2,\dots\}$ and $\N_{\sleep} = \N_0 \cup \{\sleep\}$ with $0<\sleep<1<2<\cdots$.
Define
$|\sleep|=1$, $\llbracket \sleep \rrbracket = 0$,
and
$|n|= \llbracket n \rrbracket = n$ for $n\in \N_0$.
Also define $\sleep+1=2$ and
\[
n \cdot \sleep =
\begin{cases}
n ,& n \ge 2, \\
\sleep ,& n = 1, \\
\text{undefined} ,& n=0.
\end{cases}
\]

The state of the ARW at time $t\geqslant 0$ is given by $\eta_t \in \Sigma = (\N_{\sleep})^{\Z^d}$, and the process evolves as follows.
For each site $x$, a Poisson clock rings at rate $(1+\lambda) \llbracket \eta_t(x) \rrbracket$.
When this clock rings, the system goes through the transition $\eta\to \tau_{x\sleep}\eta$ with probability $\frac{\lambda}{1+\lambda}$, otherwise $\eta\to \tau_{xy}\eta$ with probability $p(y-x)\frac{1}{1+\lambda}$.
These transitions only occur if $\eta(x)>0$ are given by
\[
  \tau_{xy}\eta(z) =
  \begin{cases}
    \eta(x)-1, & z=x, \\
    \eta(y)+1, & z=y, \\
    \eta(z),   & \mbox{otherwise,}
  \end{cases}
\quad
  \tau_{x\sleep}\eta(z) =
  \begin{cases}
    \eta(x) \cdot \sleep, & z=x, \\
    \eta(z),   & \mbox{otherwise}.
  \end{cases}
\]

We assume that $\eta_0(x)\in\N_0$ for all $x$ a.s., and use $\PP^\nu$ to denote the law of $(\eta_t)_{t\geqslant0}$, where $\nu$ denotes the distribution of $\eta_0$.
We say that $(\eta_t)_{t\geqslant0}$ \emph{fixates} if $\eta_t(x)$ is eventually constant for each fixed $x\in\Z^d$.

\paragraph{Site-wise representation and stabilization}

We now use $\eta$ to denote configurations in $\Sigma$ instead of a continuous-time process.
We say that site $x$ is \emph{unstable for the configuration $\eta$} if $\eta(x) \geqslant 1$. Otherwise, $x$ is said to be  \emph{stable}.
By \emph{toppling} site $x$ we mean the application of an operator $\tau_{xy}$ or $\tau_{x\sleep}$ to $\eta$.
Toppling an unstable site is \emph{legal}.
If $\eta(x)=\sleep$, toppling $x$ is not legal but is \emph{acceptable} (here we are departing from the original dynamics of the model, but this operation is useful in the proofs), and to that end we define $\sleep-1=0$ and $\sleep\cdot\sleep=\sleep$.
Legal topplings are acceptable.
If $\eta(x)=0$, toppling $x$ is not acceptable.

Let $\Ir=(\tau^{x,j})_{x\in\Z^d, j\in\N}$ be a fixed field of instructions, that is, for each $x$ and $j$, $\tau^{x,j}$ equals $\tau_{x\sleep}$ or $\tau_{xy}$ for some $y$.
Let $h\in (\N_0)^{\Z^d}$.
This field $h$ is called the \emph{odometer} and it counts how many topplings occur at each site.
The toppling operation at $x$ is defined by
\(
\Phi_x(\eta,h)= \big(\tau^{x,h(x)+1}\eta, h + \delta_x \big).
\)
Given a finite sequence $\alpha=(x_1,\dots,x_k)$, define 
\(
\Phi_\alpha = \Phi_{x_k}\circ \Phi_{x_{k-1}}\circ \cdots \circ \Phi_{x_1}.
\)
We say that $\alpha$ is a \emph{legal} or \emph{acceptable} sequence of topplings for $(\eta,h)$ if, for every $j=1,\dots,k$, $\Phi_{x_j}$ is legal or acceptable for $\Phi_{x_1,\dots,x_{j-1}}(\eta,h)$, respectively.
Given $V \subseteq \Z^d$, we say that $(\eta,h)$ is \emph{stable in $V$} if every $x\in V$ is stable for $\eta$.
We write $\alpha \subseteq V$ if $x_1,\dots,x_k \in V$.
We say that $\alpha$ \emph{stabilizes $(\eta,h)$ in $V$} if $\alpha$ is acceptable for $(\eta,h)$ and $\Phi_\alpha(\eta,h)$ is stable in $V$.
Let $m_\alpha$ be given by 
\(
m_\alpha(x)=\sum_\ell \I_{\{x_\ell=x\}}, \ \forall x \in \Z^d.
\)
We write $m_\beta  \leqslant m_\alpha$ if $m_\beta(x)  \leqslant m_\alpha(x)$ for all $x \in \Z^d$, and the same for $\eta \leqslant \tilde\eta$.
The following lemmas are proved in~\cite{Rolla15}.

\begin{lemma}
[Local abelianness]
If $\alpha$ and $\beta$ are acceptable sequences of topplings for the configuration $(\eta,h)$, such that $m_\alpha=m_\beta$, then $\Phi_\alpha(\eta,h) = \Phi_\beta(\eta,h)$.
\end{lemma}

For $V\subseteq \Z^d$, let
\[
m_{V,\eta,h}(x) = \sup\{ m_\beta(x) : \beta\subseteq V\text{ legal for } (\eta,h) \}
,
\]
$x \in \Z^d$.
Notice that 
$m_{V,\eta,h}(x) = \sup\{ m_{V',\eta,h}(x) : V' \subseteq V \text{ finite} \}$.
Let $m_{\eta,h} = m_{\Z^d,\eta,h}$.

\begin{lemma}
[Least Action Principle]
If $\alpha$ is an acceptable sequence of topplings that stabilizes $(\eta,h)$ in $V$,
then $m_{V,\eta,h} \le m_\alpha$.
\end{lemma}

\begin{lemma}
[Global abelianness]
If $\alpha$ and $\beta$ are both legal toppling sequences for $(\eta,h)$ that are contained in $V$ and stabilize $(\eta,h)$ in $V$, then $m_\alpha=m_\beta=m_{V,\eta,h}$.
In particular, $\Phi_\alpha(\eta,h)=\Phi_\beta(\eta,h)$.
\end{lemma}

\begin{lemma}
[Monotonicity]
If $V\subseteq \tilde V$ and $\eta \leqslant \tilde\eta$, then $m_{V,\eta,h}\leqslant m_{\tilde V,\tilde\eta,h}$.
\end{lemma}

\begin{definition*}
A configuration $\eta$ is said to be \emph{stabilizable} starting from odometer $h$ if $m_{\eta,h}(x)<\infty$ for every $x\in\Z^d$, and it is said to be \emph{explosive} if $m_{\eta,h}(x)=\infty$ for every $x\in\Z^d$.
The field of instructions $\Ir$ can be made explicit in the notation, for example, we say that $\eta$ is $\Ir$-stabilizable if $m_{\eta,h;\Ir}(x)<\infty$ for every $x\in\Z^d$.
The letter $h$ can be omitted in all the above notation when $h$ is identically zero.
\end{definition*}

\paragraph{Stabilization and fixation}



The previous properties are true for any fixed $\Ir$.
From now on we will take $\Ir$ to be random, and distributed as follows.
For each $x\in\Z^d$ and $j\in\N$, choose $\tau^{x,j}$ as $\tau_{xy}$ with probability $\frac{p(y-x)}{1+\lambda}$ or $\tau_{x\sleep}$ with probability $\frac{\lambda}{1+\lambda}$, independently over $x$ and $j$.
Let $\eta_0 \in \N_0^{\Z^d}$ have a spatially ergodic distribution $\nu$ with finite density $\nu(\eta_0(o))<\infty$, and be independent of $\Ir$.
To avoid extra notation, we define the field $\Ir$ on the same probability space $\PP^\nu$.

\begin{lemma}
$\PP^\nu(\text{fixation of }(\eta_t)_{t\geqslant 0})
= \PP^\nu ( \eta_0 \text{ stabilizable})
= \PP^\nu (m_{\eta_0}(o)<\infty)
= 0 \text{ or } 1$.
\end{lemma}

\paragraph{Mass-transport principle}

The mass-transport principle will be used a couple of times in the proof.
It consists in the following.
Let $f(x,y) = f(x,y;\omega)$ be a non-negative function of two points $x,y\in \Z^d$ and the randomness $\omega$.
Suppose that translations $\theta$ of $\Z^d$ define a group action on the randomness $\omega$, and that the law of $\omega$ is invariant under this action.
If $f(\theta x,\theta y; \theta \omega)=f(x,y; \omega)$ for every $x,y,\omega$, then
$$ \E \sum_y f(x,y) = \E \sum_y f(y,x), \qquad \text{for any } x\in \Z^d.$$
This identity says that on average the mass sent by $x$ equals the mass received by $x$, although ``mass'' can be any non-negative function.
See~\cite{LyonsPeres16} for further explanation, proof and applications.

\section{Stabilization}
\label{sec:stabilization}

In this section, we state a theorem equivalent to Theorem~\ref{thm:main}, briefly sketch its proof, and then give the proof in three parts: embedding the initial configuration into another one with higher density, stabilization of the embedded configuration, and finally stabilization of the original configuration.

Instead of proving Theorem~\ref{thm:main} as stated, we consider the following equivalent formulation.

\begin{theorem}
Let $d$, $\lambda$ and $p(\cdot)$ be given.
Let $\nu_1$ and $\nu_2$ be two spatially ergodic distributions for initial states on $\Z^d$, with respective densities $\zeta_1<\zeta_2$.
If the ARW system is a.s.\ fixating with initial state ${\nu_2}$, then it is also a.s.\ fixating with initial state ${\nu_1}$.
\end{theorem}

Let us give a brief sketch before moving to the proof.

The proof is algorithmic and has two stages, both stages being infinite. The idea is very simple and is related to what is sometimes called decoupling. Let $\eta_0$ and $\xi_0$ be independent and distributed as $\nu_1$ and $\nu_2$.
In the first stage, we evolve $\eta$ starting from $\eta_0$ until it gives a configuration $\eta_0' \le \xi_0$. In the second stage, we use the same set of instructions to evolve both systems.
Since the evolution of $\xi$ starting from $\xi_0$ fixates a.s., so does the evolution of $\eta$ starting from $\eta_0'$, concluding the proof. More precisely, in the first stage we force each particle in the system $\eta$ to move (by waking it up if needed) until it meets a particle of $\xi_0$; once they meet, they are paired and will not be moved until the second stage.
Even if it takes infinitely many steps to finish pairing globally, a.s.\ every particle in the system $\eta$ will eventually be paired, and the resulting odometer will be a.s.\ finite at every site (if the odometer were infinite somewhere, by ergodicity it would be infinite everywhere, so every particle in $\xi_0$ would be paired, implying $\zeta_2 \le \zeta_1$)
This yields a configuration $\eta_0' \le \xi_0$.
In the second stage, we simply evolve the system using the remaining instructions. Since they are independent of $\xi_0$, $\eta_0$, and of the instructions used in the first stage, by assumption the remaining instructions a.s.\ stabilize $\xi_0$ leaving a locally-finite odometer. By monotonicity of the final odometer with respect to the configuration, the same set of remaining instructions also stabilizes $\eta_0'$, again with a locally-finite odometer. Adding the odometer of both stages would give the final locally-finite odometer given by stabilization of $\eta_0$, except that in the first stage we have not followed the toppling rules correctly. But it still gives an upper bound due to monotonicity of the odometer with respect to waking up particles.

We now turn to the proof.
To make the argument precise we will not `move' particles as in the previous sketch, as the embedding requires an infinite number of topplings.
We instead explore the instructions and define a sequence of configurations in terms of $\eta_{0}$, $\xi_0$ and $\Ir$.
We end up concluding that a.s.\ the result of this exploration implies that $\eta_0$ is stabilizable, which in turn implies the statement of the theorem.

\paragraph{Embedding of the smaller configuration}
Without loss of generality we can assume that $\nu_1$ or $\nu_2$ is not only ergodic but also mixing, otherwise consider $\nu_3$ as i.i.d.\ Poisson with mean $\frac{\zeta_1+\zeta_2}{2}$ which is mixing, and apply the result from $\nu_2$ to $\nu_3$ and from $\nu_3$ to $\nu_1$.
From this assumption we get that the triple $(\eta_0,\xi_0,\Ir)$ is spatially ergodic, which is crucial in the argument.

Let $\eta_0$, $\xi_0$ and $\Ir$ be given, and take $h_0 \equiv 0$.
For $k=1,2,3,4,\dots$, suppose $\eta_{k-1}$ and $h_{k-1}$ have been defined.
Denote by $A_k$ the set given by
\[
A_k = \{ x : \eta_{k-1}(x)>\xi_0(x) \}
,
\]
and consider an arbitrary enumeration $$A_k=\{x_1^k,x_2^k,x_3^k,\dots\}.$$
Let
$$(\eta_k,h_k) = \lim_j \Phi_{(x_1^k,x_2^k,\dots,x_j^k)} (\eta_{k-1},h_{k-1})$$
in case $A_k$ is infinite -- in case it is finite, by ergodicity it is a.s.\ empty in which case we let $(\eta_k,h_k) = (\eta_{k-1},h_{k-1})$.
Notice that the condition $\eta_{k-1}(x)>\xi_0(x)$ is also satisfied when $\eta_{k-1}(x)=\sleep$ and $\xi_0(x)=0$, so this operation may require waking up particles.

As we go through $j=1,2,3,\dots$ in the above expression, for each $j$ the field $m_{(x_1^k,x_2^k,\dots,x_j^k)}$ is increased by one unit at $x_j^k$, so $h_k$ is well-defined and satisfies
$$h_k(x) = h_{k-1}(x) + \I_{A_k}(x).$$
The limit taken in $j$ to obtain $\eta_k$ from $(\eta_{k-1},h_{k-1})$ is also well-defined because, for each site $x$, the sequence decreases for at most one value of $j$.
In case it decreases, it may send one particle to one other site $z \ne x$.
We claim that the configuration at each site $x$ increases a finite number of times, hence the limit $\eta_k$ is a.s.\ finite and by the Abelian property it does not depend on the ordering of $A_k$.
Indeed, let $f(x,y)$ be the indicator of the event that $x\in A_k$ and toppling $x$ sends a particle to $y$.
Since $\sum_y f(x,y) = \I_{x\in A_k} \le 1$, by the mass-transport principle we have $\E \sum_y f(y,x) \le 1$, so $\sum_y f(y,x) < \infty$ a.s., proving the claim.

We now prove that, if $\lim_k h_k(o)=\infty$ with positive probability then we must have $\zeta_1 \geq \zeta_2$.

First, we argue that $\PP(h_0'(o)=\infty)=0 \text{ or } 1$, where $$h_0'(x)=\lim_k h_k(x).$$
By the Local Abelianness, $(\eta_k,h_k)$ does not depend on the enumeration of $A_k$.
In particular, $\eta_k$ and $h_k$ are determined by $\eta_{0}$, $\xi_0$ and $\Ir$ in a translation-covariant way, so the random set of sites $x$ for which $h_0'(x)=\infty$ is ergodic with respect to translations.
Moreover, the event $h_0'(o)=\infty$ a.s.\ implies the event that $h_0'(z)=\infty$ for every $z$ such that $p(z)>0$. These two facts together imply that, either $h_0'(x)=\infty$ a.s.\ for every $x$, or $h_0'(x)<\infty$ a.s.\ for every $x$, see~\cite[proof of Lemma~4]{RollaSidoravicius12}.
This proves the zero-one law.

Suppose $h_0'(o)=\infty$ with positive probability.
By the zero-one law we have $h_0'(o)=\infty$ a.s., which means that $\PP(\limsup_k \{o \in A_k\})=1$.
But if $o \in A_{k_0}$ for some $k_0$, then necessarily $\eta_{k_0-1}(o) > \xi_0(o)$ and hence, by definition of $A_k$, $\eta_{k}(o) \ge \xi_0(o)$ for all $k \ge k_0$ and therefore $\liminf_k \|\eta_k(o)\| \ge \|\xi_0(o)\|$.
On the other hand, from the mass-transport principle we have $\E \|\eta_k(o)\| = \E \| \eta_{k-1}(o) \| = \dots = \E \| \eta_0(o) \| = \zeta_1$ (to show the first identity, we let $f_k(x,y)$ be the indicator that, on step $k$, $x$ sends a particle to $y$, and let $f_k(x,x)$ be the number of particles that were present at $x$ at the beginning of stage $k$ and stayed at $x$).
By Fatou's Lemma, $\zeta_1 \ge \zeta_2$.

Since we are assuming $\zeta_1<\zeta_2$, we must have $h_0'(o)<\infty$ a.s.
Now, as we go through $k=1,2,3,\dots$, the value  of $\eta_k(o)$ can decrease only when $o\in A_k$, i.e., only when $h_k(o)$ increases.
Hence, $(\eta_k(o))_{k}$ is a.s.\ eventually non-decreasing, so it converges.
Its limit $\eta_0'(o)$ satisfies $\eta_0'(o) \le \xi_0(o)$, otherwise $o$ would be in $A_k$ for all large enough $k$ and $h_0'(o)$ would be infinite.
By translation invariance, a.s.\ $h_0'(x)<\infty$ and $\eta_0'(x) = \lim_k\eta_k(x) \le \xi_0(x)$ for every $x$.


\paragraph{Stabilization of the embedded configuration}

In the previous stage we obtained a pair $(\eta_0',h_0')$ a.s.\ satisfying $h_0'(x)<\infty$ and $\eta_0'(x) \le \xi_0(x)$ for every $x \in \Z^d$.
Let $\tilde{\Ir}$ be the set of instructions given by
\[
\tilde{\tau}^{x,j} = \tau^{x,h_0'(x)+j}, \qquad x \in \Z^d, j\in\N.
\]
that is, the field obtained by deleting the instructions used in the embedding stage described above.
Since the first $h_0'(x)$ instructions have been deleted at each site $x$, stabilizing a system with the instructions in $\tilde{\Ir}$ instead of $\Ir$ is equivalent to starting with odometer at $h_0'$ instead of $h_0 \equiv 0$.

Now note that the collection of instructions $\big( \tau^{x,j} : x\in \Z^d, j > h_0'(x) \big)$ played no role in the construction of $\eta_0'$ and $h_0'$, so they are independent of $\xi_0$ and $h_0'$.
Hence, $\tilde{\Ir}$ is an i.i.d.\ field just like $\Ir$, and it is also independent of $\xi_0$.

Therefore, $\PP\big[\xi_0 \text{ is $\tilde{\Ir}$-stabilizable}\big] = \PP\big[\xi_0 \text{ is $\Ir$-stabilizable}\big]$, and the latter equals $1$ by assumption.
Since $\eta_0' \le \xi_0$, we have $\PP\big[\eta_0' \text{ is $\tilde{\Ir}$-stabilizable}\big] \ge \PP\big[\xi_0 \text{ is $\tilde{\Ir}$-stabilizable}\big] = 1$.
This means that a.s.\ there exists $h_1'$ such that, for all finite $V \subseteq \Z^d$ and all $x\in\Z^d$,
$m_{V,\eta_0';\tilde{\Ir}}(x) \le h_1'(x) < \infty.$

\paragraph{Stabilization of the smaller configuration}

We now recall some properties from these two stages to show that $\eta_0$ is a.s.\ $\Ir$-stabilizable with 
$$m_{\eta_0;\Ir}(x) \le h_0'(x) + h_1'(x) < \infty, \quad \forall\ x \in\Z^d.$$

In the first stage, the limits $\eta_0'$ and $h_0'$, which are determined by $\eta_{0}$, $\xi_0$ and $\Ir$, almost-surely exist and satisfy $h_0'<\infty$ and $\eta_0' \le \xi_0$. Suppose this event occurs, and let $V$ be a fixed finite set.

If we start from $(\eta_0,h_0)$ and perform all topplings in $V$ as well as particle additions to $V$ (coming from $V^c$), following the same order as in the first stage, only a finite number of operations will be performed, and we end up with a state that equals $(\eta_0',h_0')$ on $V$.

By the local Abelian property, we can add the particles first and then topple the sites in $V$ as in the first stage, obtaining the same result.
This means that there is some $\bar\eta_V \ge \eta_0$ and an acceptable sequence $\alpha_V=(x_1,\dots,x_n)$ for $(\bar\eta_V,h_0)$ such that
$m_{\alpha_V} = h_0'$ on $V$ and
\[
\Phi_{\alpha_V} (\bar\eta_V,h_0) = (\eta_0',h_0') \text{ on } V.
\]

Now, in the second stage, we showed that a.s.\ there exists $h_1'(x) < \infty$ such that $m_{V',\eta_0';\tilde{\Ir}}(x) \le h_1'(x)$ for any finite $V'$.
Suppose this event occurs.

Notice that $m_{V,\eta_0',h_0';{\Ir}}(x) = m_{V,\eta_0';\tilde{\Ir}}(x)$, that is, to stabilize $\eta_0'$ in $V$ using the shifted field of instructions is the same as stabilize $\eta_0'$ in $V$ using the original field of instructions and shifted odometer.
Therefore, there exists $\beta_V=(x_{n+1},\dots,x_m)$ contained in $V$ such that
$m_{\beta_V} \le h_1'$ on $V$ and
$\Phi_{\beta_V} (\eta_0',h_0')$ is stable in $V$.

By the above identity,
\(
\Phi_{\beta_V} \circ \Phi_{\alpha_V} (\bar\eta_V,h_0) = \Phi_{\beta_V} (\eta_0',h_0')
,
\)
on $V$.
Since the latter is stable in $V$, by the Least Action Principle we have
$$m_{V,\bar\eta_V;\Ir}(x) \le m_{\alpha_V}(x) + m_{\beta_V}(x), \quad \text{ for all } x\in\Z^d.$$
Thus, by monotonicity,
\[
m_{V,\eta_0;\Ir}(x)
\le
m_{V,\bar\eta_V;\Ir}(x)
\le m_{\alpha_V}(x) + m_{\beta_V}(x)
\le h_0'(x) + h_1'(x) < \infty.
\]
We now note that the above bound does not depend on $V$, so
\[
m_{\eta_0;\Ir}(x) = \sup_{V \text{ finite}} m_{V,\eta_0;\Ir}(x)
\le h_0'(x) + h_1'(x) < \infty, \quad \forall\ x \in\Z^d,
\]
which means that $\eta_0$ is stabilizable, concluding the proof.

%
%

\section*{Acknowledgements}

We thank Ronald Dickman who introduced us to this topic and has been a source of inspiration for almost twenty years.
We thank Augusto Teixeira and Lionel Levine for fruitful and encouraging discussions, and Feng Liang who found a gap in an earlier version of the proof. We also thank Roberto Fernández for valuable help.
O.Z. thanks NYU-Shanghai for support and hospitality.
O.Z. also thanks support from ANR MALIN.

{
\small
\setstretch{1}
\bibliographystyle{bib/leoabbrv}
\bibliography{bib/leo}

\begin{thebibliography}{10}
\expandafter\ifx\csname urlstyle\endcsname\relax
  \providecommand{\doi}[1]{doi:\discretionary{}{}{}#1}\else
  \providecommand{\doi}{doi:\discretionary{}{}{}\begingroup
  \urlstyle{rm}\Url}\fi

\bibitem{AmirGurel-Gurevich10}
G.~\textsc{Amir}, O.~\textsc{Gurel-Gurevich}.
\newblock \emph{On fixation of activated random walks}.
\newblock Electron Commun Probab \textbf{15}:119--123, 2010.
\newblock \doi{10.1214/ECP.v15-1536}.

\bibitem{BakTangWiesenfeld88}
P.~\textsc{Bak}, C.~\textsc{Tang}, K.~\textsc{Wiesenfeld}.
\newblock \emph{Self-organized criticality}.
\newblock Phys Rev A \textbf{38}:364--374, 1988.
\newblock \doi{10.1103/PhysRevA.38.364}.

\bibitem{BasuGangulyHoffman18}
R.~\textsc{Basu}, S.~\textsc{Ganguly}, C.~\textsc{Hoffman}.
\newblock \emph{Non-fixation for conservative stochastic dynamics on the line}.
\newblock Comm Math Phys \textbf{358}:1151--1185, 2018.
\newblock \doi{10.1007/s00220-017-3059-7}.

\bibitem{CabezasRollaSidoravicius14}
M.~\textsc{Cabezas}, L.~T. \textsc{Rolla}, V.~\textsc{Sidoravicius}.
\newblock \emph{Non-equilibrium phase transitions: Activated random walks at
  criticality}.
\newblock J Stat Phys \textbf{155}:1112--1125, 2014.
\newblock \doi{10.1007/s10955-013-0909-3}.

\bibitem{CabezasRollaSidoravicius18}
---{}---{}---.
\newblock \emph{Recurrence and density decay for diffusion-limited annihilating
  systems}.
\newblock Probab Theory Relat Fields \textbf{170}:587--615, 2018.
\newblock \doi{10.1007/s00440-017-0763-3}.

\bibitem{DickmanMunozVespignaniZapperi00}
R.~\textsc{Dickman}, M.~A. \textsc{Mu{\~n}oz}, A.~\textsc{Vespignani},
  S.~\textsc{Zapperi}.
\newblock \emph{Paths to self-organized criticality}.
\newblock Braz J Phys \textbf{30}:27, 2000.
\newblock \doi{10.1590/S0103-97332000000100004}.

\bibitem{DickmanRollaSidoravicius10}
R.~\textsc{Dickman}, L.~T. \textsc{Rolla}, V.~\textsc{Sidoravicius}.
\newblock \emph{Activated random walkers: Facts, conjectures and challenges}.
\newblock J Stat Phys \textbf{138}:126--142, 2010.
\newblock \doi{10.1007/s10955-009-9918-7}.

\bibitem{FeyLevineWilson10a}
A.~\textsc{Fey}, L.~\textsc{Levine}, D.~B. \textsc{Wilson}.
\newblock \emph{Driving sandpiles to criticality and beyond}.
\newblock Phys Rev Lett \textbf{104}:145703, 2010.
\newblock \doi{10.1103/PhysRevLett.104.145703}.

\bibitem{FeyMeester15}
A.~\textsc{Fey}, R.~\textsc{Meester}.
\newblock \emph{Critical densities in sandpile models with quenched or annealed
  disorder}.
\newblock Markov Process Related Fields \textbf{21}:57--83, 2015.
\newblock \href{http://arxiv.org/abs/1211.4760}{arXiv}.

\bibitem{Fey-denBoerRedig05}
A.~\textsc{{Fey-den Boer}}, F.~\textsc{Redig}.
\newblock \emph{Organized versus self-organized criticality in the {Abelian}
  sandpile model}.
\newblock Markov Process Related Fields \textbf{11}:425--442, 2005.
\newblock \href{http://arxiv.org/abs/math-ph/0510060}{arXiv}.

\bibitem{JoJeong10}
H.-H. \textsc{Jo}, H.-C. \textsc{Jeong}.
\newblock \emph{Comment on ``driving sandpiles to criticality and beyond''}.
\newblock Phys Rev Lett \textbf{105}:019601, 2010.
\newblock \doi{10.1103/PhysRevLett.105.019601}.

\bibitem{Levine15}
L.~\textsc{Levine}.
\newblock \emph{Threshold state and a conjecture of {P}oghosyan, {P}oghosyan,
  {P}riezzhev and {R}uelle}.
\newblock Comm Math Phys \textbf{335}:1003--1017, 2015.
\newblock \doi{10.1007/s00220-014-2216-5}.

\bibitem{LyonsPeres16}
R.~\textsc{Lyons}, Y.~\textsc{Peres}.
\newblock \emph{Probability on trees and networks}, vol.~42 of \emph{Cambridge
  Series in Statistical and Probabilistic Mathematics}.
\newblock Cambridge University Press, New York, 2016.
\newblock \doi{10.1017/9781316672815}.
\newblock \href{http://mypage.iu.edu/\string~rdlyons}{free download}.

\bibitem{PoghosyanPoghosyanPriezzhevRuelle11}
S.~S. \textsc{Poghosyan}, V.~S. \textsc{Poghosyan}, V.~B. \textsc{Priezzhev},
  P.~\textsc{Ruelle}.
\newblock \emph{Numerical study of the correspondence between the dissipative
  and fixed-energy abelian sandpile models}.
\newblock Phys Rev E \textbf{84}:066119, 2011.
\newblock \doi{10.1103/PhysRevE.84.066119}.

\bibitem{Rolla15}
L.~T. \textsc{Rolla}.
\newblock \emph{Activated random walks}, 2015.
\newblock Preprint. \href{http://arxiv.org/abs/1507.04341}{arXiv:1507.04341}.

\bibitem{RollaSidoravicius12}
L.~T. \textsc{Rolla}, V.~\textsc{Sidoravicius}.
\newblock \emph{Absorbing-state phase transition for driven-dissipative
  stochastic dynamics on {$Z$}}.
\newblock Invent Math \textbf{188}:127--150, 2012.
\newblock \doi{10.1007/s00222-011-0344-5}.

\bibitem{RollaTournier18}
L.~T. \textsc{Rolla}, L.~\textsc{Tournier}.
\newblock \emph{Non-fixation for biased activated random walks}.
\newblock Ann Inst H Poincar{\'e} Probab Statist \textbf{54}:938--951, 2018.
\newblock \doi{10.1214/17-AIHP827}.

\bibitem{Shellef10}
E.~\textsc{Shellef}.
\newblock \emph{Nonfixation for activated random walks}.
\newblock ALEA Lat Am J Probab Math Stat \textbf{7}:137--149, 2010.
\newblock \href{http://alea.impa.br/articles/v7/07-07.pdf}{pdf}.

\bibitem{SidoraviciusTeixeira17}
V.~\textsc{Sidoravicius}, A.~\textsc{Teixeira}.
\newblock \emph{Absorbing-state transition for stochastic sandpiles and
  activated random walks}.
\newblock Electron J Probab \textbf{22}:33, 2017.
\newblock \doi{10.1214/17-EJP50}.

\bibitem{StaufferTaggi18}
A.~\textsc{Stauffer}, L.~\textsc{Taggi}.
\newblock \emph{Critical density of activated random walks on transitive
  graphs}.
\newblock Ann Probab \textbf{to appear}, 2018.
\newblock \href{http://arxiv.org/abs/1512.02397}{arXiv:1512.02397}.

\bibitem{Taggi16}
L.~\textsc{Taggi}.
\newblock \emph{Absorbing-state phase transition in biased activated random
  walk}.
\newblock Electron J Probab \textbf{21}:13, 2016.
\newblock \doi{10.1214/16-EJP4275}.

\end{thebibliography}
\par
}

\end{document}